\documentclass[12pt, reqno]{article}
\usepackage{geometry}
\usepackage{graphicx}
\usepackage{amssymb,latexsym, amsmath}
\usepackage{amsfonts,amsthm}
\usepackage{amscd}
\usepackage{mathrsfs}
\usepackage{cite,upref}
\usepackage{eucal}
\usepackage{epstopdf}
\usepackage{amsgen}
\usepackage{xspace}
\usepackage{verbatim}

\newcommand{\lf}{L_{\ms{F}}^2}

\newcommand{\gd}{\Delta}

\newcommand{\gw}{\Omega}

\newcommand{\mb}{\mathbb}
\newcommand{\ms}{\mathscr}

\newcommand{\om}{\omega}

\newcommand{\nb}{\nabla}

\newcommand{\Tu}{\Tilde{u}}

\newcommand{\tbf}{\textbf}
\newcommand{\tup}{\textup}

\newcommand{\ragv}{\rangle}

\newcommand{\inptv}[1]{\lag #1 \ragv}
\newcommand{\ef}[1]{E_{\ms{F}_{#1}}}

\newcommand{\hn}[1]{\| #1 \|_H}
\newcommand{\vn}[1]{\| #1 \|_V}

\newcommand{\eh}{\emph}

\newcommand{\beq}{\begin{equation}}
\newcommand{\eeq}{\end{equation}}
\newcommand{\bea}{\begin{align}}
\newcommand{\eea}{\end{align}}
\newcommand{\bec}{\begin{center}}
\newcommand{\eec}{\end{center}}
\newcommand{\bthm}{\begin{theorem}}
\newcommand{\ethm}{\end{theorem}}
\newcommand{\bpr}{\begin{proof}}
\newcommand{\epr}{\end{proof}}
\newcommand{\bcl}{\begin{corollary}}
\newcommand{\ecl}{\end{corollary}}
\newcommand{\bpn}{\begin{proposition}}
\newcommand{\epn}{\end{proposition}}
\newcommand{\bre}{\begin{remark}}
\newcommand{\ere}{\end{remark}}
\newcommand{\bdf}{\begin{definition}}
\newcommand{\edf}{\end{definition}}
\newcommand{\bss}{\begin{align*}}
\newcommand{\ess}{\end{align*}}
\newcommand{\bst}{\begin{split}}
\newcommand{\est}{\end{split}}
\newcommand{\lag}{\langle}
\newcommand{\rag}{\rangle}
\newcommand{\bl}{\label}

\numberwithin{equation}{section}

\theoremstyle{plain}
\newtheorem{theorem}{Theorem}[section]
\newtheorem{lemma}[theorem]{Lemma}
\newtheorem{corollary}[theorem]{Corollary}
\newtheorem{proposition}{Proposition}
\theoremstyle{definition}
\newtheorem{definition}{Definition}
\theoremstyle{remark}
\newtheorem{remark}{Remark}

\DeclareMathOperator*{\esssup}{ess\,sup}

\begin{document}

\title{2D Backward
Stochastic Navier-Stokes Equations  with Nonlinear Forcing~\thanks{Supported by NSFC Grant \#10325101, by Basic Research Program of China (973
Program)  Grant \# 2007CB814904, by the Science Foundation of the Ministry of Education of China Grant \#200900071110001, and by WCU (World
Class University) Program through the Korea Science and Engineering Foundation funded by the Ministry of Education, Science and Technology
(R31-2009-000-20007). Part of this study was reported by the second author at the third International Symposium of Backward Stochastic
Differential Equations and their Applications, held in 2002 in Weihai, Shandong Province, China.  The third author thanks the School of
Mathematical Sciences, Fudan University for the hospitality during his visit in the winter of 2009.}}

\author{Jinniao Qiu~\thanks{Department of Finance and Control Sciences, School of Mathematical Sciences, Fudan University, Shanghai 200433, China.
\textit{E-mail}: \texttt{071018032@fudan.edu.cn} (Jinniao Qiu), \texttt{sjtang@fudan.edu.cn} (Shanjian Tang).}\ , \quad Shanjian
Tang~${}^\dag$\thanks{Graduate Department of Financial Engineering, Ajou University, San 5, Woncheon-dong, Yeongtong-gu, Suwon, 443-749,
Korea.}\ , \quad and Yuncheng You~\thanks{Department of Mathematics and Statistics, University of South Florida, Tampa, FL 33620, USA.
\textit{E-mail}:
 \texttt{you@math.usf.edu} (Yuncheng You).}}

\maketitle

\begin{abstract}
The paper is concerned with the existence and uniqueness of a strong solution to  a two-dimensional backward stochastic Navier-Stokes equation
with nonlinear forcing, driven by a Brownian motion.  We use the spectral approximation and  the truncation and variational techniques. The
methodology features an interactive analysis on basis of the regularity of the deterministic Navier-Stokes dynamics and the stochastic
properties of the It\^o-type diffusion processes.
\end{abstract}



\section{\textbf{Introduction}}

Let $(\gw, \ms{F}, \{\ms{F}_t\}, P)$ be a complete, filtrated probability space, on which defined is a standard 1-dimensional Brownian motion
$\{W_t\}_{t \geq 0}$, whose natural augmented filtration is denoted by $\{\ms{F}_t, t\in [0,T]\}$. We denote by $\mathcal {P}$ the
$\sigma$-Algebra of the predictable sets on $\Omega\times[0,T]$ associated with $\{\ms{F}_t\}_{t\geq 0}$. The expectation will be exclusively
denoted by $E$ and the conditional expectation on $\ms{F}_s$ will be denoted by $E_{\ms{F}_s}$. We use \emph{a.s.} to denote that an equality
or inequality holds almost surely with respect to the probability measure $P$.

The theory of backward stochastic differential equations (BSDEs) has received an extensive studies in the last two decades in connection with
a wide range of applications as in stochastic control theory, econometrics, mathematical finance, and nonlinear partial differential
equations. See \cite{DE92, DT, DL92, KM97, MPY94,PP90, YZ99} for details. In Tang~\cite{Tang}, a very general system of backward stochastic
partial differential equations (BSPDEs) are studied. However, they are semi-linear, and thus exclude the nonlinearity of the well-known
Navier-Stokes operator. In this paper, we concentrate our attentions to study the backward stochastic Navier-Stokes equation (BSNSE).

The standard deterministic Navier-Stokes equation describing the velocity field of an incompressible, viscous fluid motion in a domain of
$\mb{R}^d$ ($d = 2$ or 3) takes the following form:
\begin{equation}\label{nse}
  \left\{\begin{array}{l}
    \begin{split}
      \partial_t u - \nu \gd u + (u\cdot \nb)u + \nb p + f=0, \;\; t \geq 0;
    \end{split}\\
    \begin{split}
     \nb \cdot u = 0, \quad u(0) = u_0,
    \end{split}
  \end{array}\right.
\end{equation}
where $u = u( t,x)$ represents the $d$-dimensional velocity field of a fluid, $p = p( t,x)$ is the pressure, $\nu\in(0,\infty)$ is the
viscosity coefficient, and $f = f(t,x)$ stands for the external force. Let $(u, p)$ solve the equation \eqref{nse}. By reversing the time and
defining
$$
    \Tu (t, x) = - u (T-t,x), \quad \Tilde{p} (t, x) = p (T-t,x), \quad \textup{for} \; t \leq T,
$$
then $(\Tu, \Tilde{p})$ satisfies the following \emph{backward} Navier-Stokes equation:
\begin{equation}\label{bnse}
  \left\{\begin{array}{l}
          \partial_t \Tu + \nu \gd \Tu + (\Tu\cdot \nb)\Tu + \nb \Tilde{p} +f=0 , \;\; t \leq T;\\
         \nb \cdot \Tu = 0, \quad \Tu(T) = \Tu_0.
    \end{array}\right.
\end{equation}
Note that the time-reversing makes the original initial value problem of \eqref{nse} become a terminal value problem of \eqref{bnse}.

We shall study the following two-dimensional backward stochastic Navier-Stokes equations (briefly 2D BSNSE) in $\mathbb{R}^2$ with a spatially
periodic condition and a given terminal condition at time $T > 0$:
\begin{equation}\label{oreq}
  \left\{
    \begin{split}
          d u(t,x) +\{\nu \gd u(t,x) +
    &
        (u \cdot \nb) u(t,x)+(\sigma\cdot \nb) Z(t,x) + \nb p(t,x) \}\, dt \\
        = - f (t,x, u, Z)\, dt +& Z(t,x)\, dW_t, \quad (t,x) \in [0, T)\times \mathbb{R}^2; \\
     \tup{div} \, u(t,x)  =\ & 0, \quad (t, x) \in [0,T)\times \mathbb{R}^2; \\
     u(t, x+ae_i) =\ & u(t,x), \quad (t,x) \in [0, T)\times \mathbb{R}^2,  \;\; i=1,2; \\
     u(T, x) = \ &\xi (x), \; \; (x, \om) \in \mathbb{R}^2.
    \end{split}\right.
\end{equation}
Here $\sigma$ is a measure of ``correlation" between the Laplace and the Brownian motion, $\{e_1,e_2\}$ is the canonical basis of
$\mathbb{R}^2$, $a>0$ is the period in the $i$th direction, $u = (u_1 , u_2 )$ is the random two-dimensional velocity field of a fluid in
$\mathbb{R}^2$,
 $f$ represents the external forces which allow for feedback involving the velocity field $u$ and  the stochastic process $Z$ and may be inhomogeneous in
time. The terminal status of the velocity field is a known random field
$\xi$ on the underlying probability space. For notational convenience,
however, the variable $\om \in \gw$ in various functions and solutions will
often be omitted.

It is worth noting that, though sharing the same name, our BSNSE essentially differs from that of Sundar and Yin~\cite{SY09}  since the sign of the
viscous term ``$\nu \gd u$"  differs. Furthermore, we allow the external force $f$ to depend on both unknown fields $u$ and $Z$ in a nonlinear
way, and the drift term to depend on the gradient of the second unknown field $Z$.

In \cite{MR04,MR05}, Cauchy problems for the (forward) stochastic Navier-Stokes equations in $\mathbb{R}^d$ driven by a random nonlinear force
and a  white noise are studied and the existence and uniqueness of a global martingale solution have been proved. As a motivation BSNSEs
emerge in regard to inverse problems to determine the stochastic noise coefficients from the terminal velocity field as observed. In
\cite{gI06,xZ09}, a stochastic representation in terms of Lagrangian paths for the backward incompressible Navier-Stokes equations without
forcing is shown and used to prove the local existence of solutions in weighted Sobolev spaces and the global existence results in two
dimensions or with a large viscosity. In~\cite{SY07,SY09}, the existence and uniqueness of adapted solutions are given to the backward
stochastic Lorenz equations and to the backward stochastic Navier-Stokes equations~\eqref{oreq} in a bounded domain with $\sigma\equiv 0$, $\nu<0$ and
the external force $f(t,y,z)\equiv f(t)$.

The rest of the paper is organized as follows. In Section 2, we introduce some notations, assumptions, and preliminary lemmas, and state the
main result (see Theorem 2.1). In Section 3, we consider the spectral approximations and give relevant estimates. In Section 4, we prove the
existence of an adapted solution to the projected finite dimensional systems for our 2D BSNSE. Finally, in Section 5, we give the proof of
Theorem 2.1.

\section{\textbf{Preliminaries and the main results}}

Let $G=(0,a)^2$ be the rectangular of the period. For any nonnegative integer $m$, we denote by $H^m(G)$ the Sobolev space of functions which
are in $L^2(G)$, together with all their derivatives of orders up to $m$ and by $H_{pe}^m(G)$ the space of functions which belong to
$H_{loc}^m(\mathbb{R}^2)$ (i.e., $u|_\mathcal {O}\in H^m(\mathcal {O})$ for every open bounded set $\mathcal{O}$) and which are periodic with
period $G$. $H^m_{pe}(G)$ is a Hilbert space for the scalar product and the norm
$$
  (u,v)_m=\sum_{[\alpha]\leq m}\int_G D^\alpha u(x)D^\alpha v(x)\ dx,\quad
  |u|_m=\{(u,u)_m\}^{1/2},
$$
where $\alpha=(\alpha_1,\alpha_2)$,
$[\alpha]=\alpha_1+\alpha_2$, and
$$
   D^\alpha=D_1^{\alpha_1} D_2^{\alpha_2}=\frac{\partial ^{[\alpha]}}{\partial
   x_1^{\alpha_1}\partial x_2^{\alpha_2}}
$$
with $\alpha_1$ and $\alpha_2$ being two nonnegative integers.
The elements of $H^m_{pe}(G)$ are characterized by their Fourier series
expansion:
\begin{equation}\label{Fourier expa}
  H^m_{pe}=\{ u:\ u=\sum_{k\in \mathbb{Z}^2}c_ke^{2i\pi k\cdot x/a},\ \bar{c_k}=c_{-k},
    \sum_{k\in \mathbb{Z}^2}|k|^{2m}|c_k|^2<\infty \},
\end{equation}
and the norm $|u|_m$ is equivalent to the norm $\left\{\sum_{k\in \mathbb{Z}^2}
(1+|k|^{2m})|c_k|^2 \right\}^{1/2}$. Set a Hilbert subspace of $H^m_{pe}(G)$:
\begin{equation}\label{Fourier expa Hom}
  \dot{H}^m_{pe}(G)=\{u\in H^m_{pe}(G):
   \textrm{ in its Fourier expansion  \eqref{Fourier expa}, }c_0 =0\},
\end{equation}
with the norm $|u|_{m,0}=\{ \sum_{k\in
\mathbb{Z}^2}|k|^{2m}|c_k|^2\}^{1/2}$.
 Actually, through \eqref{Fourier expa} and \eqref{Fourier
expa Hom}, we can define $H^m_{pe}(G)$ and $\dot{H}^m_{pe}(G)$ for
arbitrary $m\in \mathbb{R}$. Moreover, $\dot{H}^m_{pe}(G)$ and
$\dot{H}^{-m}_{pe}(G)$ are in duality for all $m\in \mathbb{R}$.

As in the framework of treating the deterministic Navier-Stokes equations
(c.f. \cite{rT88,rT95,SY02}), we set up three phase spaces of functions of
the spatial variable $x \in G$ as follows:
\begin{align*}
    H &= \{\varphi \in \dot{H}^0_{pe} (G) \times \dot{H}^0_{pe} (G): \tup{div}\,
    \varphi = 0 \textrm{ in } \mathbb{R}^2 \}, \\
    V &= \{\varphi \in \dot{H}_{pe}^1 (G) \times \dot{H}_{pe}^1 (G): \tup{div} \,
    \varphi = 0\textrm{ in } \mathbb{R}^2\},\\
    \mathscr{V}&=V\bigcap C^\infty(\mathbb{R}^2)\times C^\infty(\mathbb{R}^2),
\end{align*}
where $C^\infty(\mathbb{R}^2)$ is the set of smooth functions on
$\mathbb{R}^2$.
Then both  $H$ and $V$ are Hilbert spaces equipped with the respective scalar product and norm
\begin{equation}
  \begin{split}
        \langle \phi,\ \varphi\rangle_{H}:
    &
        =  \sum_{i=1}^{2} (\phi^i,\ \varphi^i)_0  ,\ \|\phi\|_H:= \{ \langle \phi,\ \phi\rangle_H \}^{1/2},\ \phi,\varphi\in H;
    \\
        \langle \phi,\ \varphi\rangle_{V}:
    &
        =  \sum_{i,j=1}^{2} (D_j^1 \phi^i,\ D_j^1 \varphi^i)_0  ,\ \|\phi\|_V:= \{ \langle \phi,\ \phi\rangle_V \}^{1/2},\
        \phi,\varphi\in V.
  \end{split}
\end{equation}
For simplicity, we denote $\|\cdot\|$ and $\langle\cdot,\ \cdot\rangle_H$ by $\|\cdot\|$ and $\langle\cdot,\ \cdot\rangle$ respectively.
 The dual product of $\psi \in V'$ and $\varphi \in V$
will be denoted by $\lag \psi , \varphi \rag_{V',V}$ and it follows that
$$\lag
\phi_1 , \phi_2 \rag_{V',V}=\lag \phi_1, \phi_2 \rag \quad \phi_1\in V,\phi_2\in
H.$$
For a little notational abuse we still denote by $\lag\cdot,\ \cdot\rag$ the dual
product $\lag\cdot,\ \cdot\rag_{V',V}$.
We shall use $| \cdot |$ to denote the absolute value or the Euclidean
norm of $\mathbb{R}^2$. The set of all positive integers will be denoted by
$\mathbb{Z}^+$ or $\mathbb{N}$. The Lebesgue measure of the domain $G$ will be denoted by
$|G|$.
 Define
$\dot{\mathbb{H}}^m:=\dot{H}_{pe}^m(G)\times \dot{H}_{pe}^m(G) $.

For any finite dimensional vector space $F$ and $a,b\in F$, we denote by $a\cdot b$ the scalar product on $F$.
Throughout this paper, we assume that the external force term $f(t,u,Z)$
and the terminal value term $\xi$ are $\dot{\mathbb{H}}^0 $-valued and $
\dot{\mathbb{H}}^1$-valued respectively, so the solution pair $(u,Z)$ of
\eqref{oreq} must be $\dot{\mathbb{H}}^0\times\dot{\mathbb{H}}^0$-valued.
By applying the projection $\mathbb{P}: \dot{\mathbb{H}}^0 \to H$ (see
\cite{rT95}), since
$$
    H^{\bot} = \{\psi \in \dot{\mathbb{H}}^0: \psi = \nb q  \; \;
    \tup{for some} \;\; q \in H^1_{pe} (G)\},
$$
we can formulate the above terminal value problem of the 2D BSNSE \eqref{oreq} into the following problem
to solve the backward stochastic evolutionary Navier-Stokes equation,
\begin{equation}\label{pb}
  \left\{\begin{array}{l}
    \begin{split}
      -d u(t) =& \{-\nu Au(t) +B(u(t))+JZ(t) + f(t, u(t), Z(t))\} dt\\
                & - Z(t) dW_t, \quad t \in [0, T),
      \\
        u(T) =& \xi,
    \end{split}
  \end{array}\right.
\end{equation}
where
\begin{equation*}
    \begin{split}
    \Pi(u, v) = \mathbb{P} ((u \cdot \nb) v): V \times V \to V^\prime ,  \quad
    B(u) = \Pi(u, u): V \to V^\prime,\\
  JZ=\mathbb{P}\left((\sigma\cdot\nb) Z\right),\quad \sigma(t,x)=(\sigma^{1}(t,x),\ \sigma^2(t,x)),\quad \left((\sigma\cdot\nb) Z\right)^i:=\sum_{j=1}^{2}\sigma^{j}Z^{i}_{x_j},
    \end{split}
\end{equation*}
 and
$$ A \varphi = \mathbb{P} (- \gd \varphi)=-\Delta\varphi,
$$
whose domain is $D(A) = \dot{\mathbb{H}}^2 \bigcap H$ and by Poincar\'{e}
inequality we can show that $V = D(A^{1/2})$. Accordingly we shall
  adopt the equivalent norm $\| \varphi \|_V =
\|\nb \varphi \|=\|A^{1/2}\varphi\|$.
Then all $H, V$, and $D(A)$
(with the graph norm $\|\cdot\|_{D(A)}$) are separable Hilbert spaces. With a little
notational abuse we still use $f$ and $Z$ for the projections
$\mathbb{P}(f)$ and $\mathbb{P}(Z)$, respectively,

Note that the above Stokes operator $A: D(A) \to H$ is positive definite, self-adjoint, and linear, and its resolvent is compact.
Therefore, all the eigenvalues of $A$ can be ordered into the increasing sequence $\{\lambda_i\}_{i=1}^{\infty}$. The corresponding
eigenfunctions $\{e_i\}_{i = 1}^{\infty}$ form a complete orthonormal basis for the space $H$, which is also a complete orthogonal (but not
orthonormal) basis of the space $V$. With the identification $H = H'$ by the Riesz mapping, one has the triplet structure of compact
(consequently continuous) embedding,
\begin{equation} \bl{trpt}
    V \subset H \subset V'.
\end{equation}

In what follows, $C>0$ is a constant which may vary from line to line and we denote by $C(a_1,a_2,\cdots)$ or $C_{a_1,a_2,\cdots}$  a constant to depend on the parameters  $a_1,a_2,\cdots$.

We make the following three assumptions.

\textbf{Assumption} (A1).  The $H$-valued mapping $f$ is defined
on $\Omega\times[0, T] \times V\times H $ and for any $(u,z)\in V\times H$,
$f(\cdot,u,z)$ is a predictable and $H$-valued process.
Moreover, there exist a nonnegative constant $\beta$ and a nonnegative adapted process
$g\in L^\infty(\Omega,L^1([0,T]))$ such that the following conditions hold
for all $v,v_1,v_2\in V$, $\phi,\phi_1,\phi_2\in H$ and
$(\omega,t)\in\Omega\times [0,T]$:

(1). the map $s\mapsto \lag f(t,v_1+sv_2,\phi),\  v  \rag  $ is continuous on
$\mathbb{R};$

(2).
\begin{equation*}
  \begin{split}
    &\lag f(t,v_1,\phi_1)-f(t,v_2,\phi_2),\  v_1-v_2 \rag \\
    \leq&\
        \rho(v_2)\left(\|v_1-v_2\|^2+\|v_1-v_2\|(\|\phi_1-\phi_2\|+\|v_1-v_2\|_V)\right);
  \end{split}
\end{equation*}
where $\rho :V\rightarrow (0,+\infty)$ is measurable and locally
bounded;

(3). $$ \lag f(t,v,\phi),\ v   \rag \leq g(t)+\epsilon
\|\phi\|^2+\varrho(\epsilon)\|v\|^2+\beta\|v\|_V\|v\| ,$$ where $\varrho :
(0,1]\rightarrow \mathbb{R}^+$ is  continuous and decreasing;

(4). $$ \|f(t,v,\phi)\|^2\leq \left(g(t)+\beta(\|v\|_V^2+\|\phi\|^2)\right)\rho_1(v)
,$$ where $\rho_1 :H\rightarrow (0,+\infty)$ is  measurable  and
locally bounded.

\begin{remark}\label{rmk Lip wrt Z}
  In fact, (1) and (2) of Assumption (A1) implies $f(t,x,u,z)$ is locally
  Lipschitz continuous with respect to $z$ in the following sense:
  $$\|f(t,u,z)-f(t,u,Z)\|_{V'}\leq C(\|u\|_V) \|z-Z\|,$$
  for all $(\omega,t)\in \Omega\times [0,T],u\in V$ and $z,Z\in H$.
  Actually, for any $\phi\in V-\{0\},\epsilon\in \mathbb{R}^+,$
   $$   \lag f(t,u+\epsilon \phi,z)-f(t,u,Z),\ {\phi} \rag
      \leq C(\|u\|_V){\|\phi\|_V} \{\epsilon \|\phi\|+\|z-Z\|\}.$$
Letting $\epsilon \downarrow 0$, from the arbitrariness of $\phi$ we
conclude that the local Lipschitz continuity holds.
\end{remark}

 \textbf{Assumption} (A2). The function $\sigma^{j}$ defined on
$\Omega\times [0,T]$ is real-valued $\mathcal{P}$-measurable such that $
|\sigma^{j}|\leq \Lambda,$ almost surely for $j=1,2$ and all $t\in
[0,T],$ for some $\Lambda\in (0,\infty)$.

\textbf{Assumption} (A3) (super-parabolicity). There exist two constants
$\lambda>0$ and $\bar{\lambda}>1$ such that
$$ 2\nu |\xi|^2-\bar{\lambda}^2 \left(\sigma(t)\cdot \xi\right)^2
\geq 2\lambda |\xi|^2 $$
holds almost surely for all
$(t,\xi)\in [0,T] \times \mathbb{R}^2$.

Note that, in Assumptions (A2) and (A3), our $\sigma$ is defined independent of the spatial variable $x$.

For  Banach space $B$ and $p>1$, define
$$
    \ms{L}_{\ms{F}}^p (0, T; B) := \{\phi \in L^p (\gw \times [0,T];
    B) \, | \, \{\phi (\cdot, t)\}_{0 \leq t \leq T}
     \, \textup{is a predictable process} \}.
$$
Define
\begin{align*}
     M [0,T] := L_{\ms{F}}^2 (\gw; C([0, T]; H)) \cap
    \ms{L}_{\ms{F}}^2 (0, T; V)
\end{align*}
 and
$$
    \ms{M} := \ms{M} [0, T]: = M [0, T] \times \ms{L}_{\ms{F}}^2 (0, T; H)
$$
 equipped with the norm
$$
    \|(u, Z)\|_{\ms{M}} = \left\{E \left[\sup_{t \in [0, T]} \| u
    (t) \|^2\right] + E \left[\int_0^T \| u(t) \|_V^2 \, dt\right] + E \left[\int_0^T \|
    Z(t) \|^2 dt \right]\right\}^{1/2}.
$$

Throughout the paper, define
\begin{equation}\label{lemma2.2 0}
 \Phi(t,\phi,\varphi):=-\nu A\phi+B(\phi)+J\varphi+f(t,\phi,\varphi),
  \ (\phi,\varphi)\in V\times H.
\end{equation}

\begin{definition}(weak solutions) For $\xi\in L_{\ms{F}_T}^\infty (\gw; H)$ given, we say that  $(u,Z)\in \ms{M}$ is
a weak solution to \eqref{pb} if for any $\varphi\in
\mathscr{V}$,  there holds almost surely
    \begin{equation}\label{def1 eq}
        \begin{split}
           \lag u(t),\ \varphi\rag=&
            \ \lag \xi,\ \varphi\rag
            +\int_t^T \lag \Phi(s,u(s),Z(s)),\ \varphi \rag\  ds\\
             &\  -\int_t^T \lag Z(s),\ \varphi \rag\  dW_s,\ \forall t\in [0,T].
        \end{split}
    \end{equation}
\end{definition}

\begin{definition}\label{def strong solution}(strong solutions)
For $\xi\in L_{\ms{F}_T}^\infty (\gw; V) $ given, we say that $(u,Z)$ is
  a strong solution to \eqref{pb} if $(u,Z)$ is a weak solution and
  $$(u,Z)\in \left(L_{\ms{F}}^2
   (\gw; C([0, T]; V)) \cap
    \ms{L}_{\ms{F}}^2 (0, T; D(A))\right)\times \ms{L}_{\ms{F}}^2 (0, T; V).$$
\end{definition}
\begin{remark}\label{rmk def}
  If we have verified that $(u,Z)\in (L_{\ms{F}}^2
   (\gw; C([0, T]; V)) \cap
    \ms{L}_{\ms{F}}^2 (0, T; D(A)))\times \ms{L}_{\ms{F}}^2 (0, T; V)$ and
    that
$$
u(t)=\xi+\int_t^T \Phi(s,u(s),Z(s))\ ds-\int_t^T  Z(s)\ dW_s\quad a.s.\textrm{
in }H,
$$
By the stochastic Fubini theorem (see, \cite[Theorem 4.18]{Prato-Zabczyk-1992}), we can check that $(u,Z)$ is  a strong solution to \eqref{pb}.
\end{remark}

The main result of the paper is stated in the following theorem.

\bthm \label{thm main}
    Under Assumptions \tup{(A1)-(A3)}, for any $\xi \in
    L_{\ms{F}_T}^\infty (\gw; V)$, the 2D BSNSE problem \eqref{pb} admits a
    unique strong solution such that
    \begin{equation}\label{main thm estimate}
      \begin{split}
       &
            \esssup_{(\omega,s)\in\Omega\times [0,T]}\|u(s)\|_V^2+ E\left[\int_0^T \|u(s)\|_{D(A)}^2 ds+
            \int_0^T \|Z(s)\|_V^2 ds  \right] \\
       &\leq
            C\left\{\|g\|_{L^\infty(\Omega,L^1([0,T]))}+\|\xi\|_{L_{\ms{F}_T}^\infty (\gw; V)}^2 \right\},
      \end{split}
    \end{equation}
    where $C$ is a constant depending on
    $\nu,\lambda,\bar{\lambda},\beta, \varrho,\rho_1 \textrm{ and }
     T$.
\ethm

For the trilinear mapping
$$
    b(u, v, w) := \inptv{\Pi(u, v),\  w} = \sum_{i= 1}^2 \sum_{j = 1}^2 \int_G u_i \, \frac{\partial
    v_j}{\partial x_i} \, w_j \, dx,\  u,v,w\in H,
$$
we have the following instrumental regularity properties.

\begin{lemma} \bl{L:bprop}
    The following properties hold for any two-dimensional bounded domain $G$, where $C_G$
    is used to denote different constants only depending on $G$.
\beq \bl{buvw}
    \begin{split}
    |b(u, v, w)| &\leq 2^{1/2} \|u\|_H^{1/2} \| u\|_V^{1/2} \| v \|_V
    \| w \|_H^{1/2} \| w \|_V^{1/2}, \quad \; \; u, v, w \in V, \\[2pt]
    |b(u, v, w)| &\leq C_G \|u\|_H^{1/2} \| Au\|_H^{1/2} \| v \|_V
    \| w \|_H, \quad \quad \quad \quad \, u \in D(A), v \in V, w \in H, \\[2pt]
    |b(u, v, w)| &\leq C_G \|u\|_H^{1/2} \| u\|_V^{1/2} \| v
    \|_V^{1/2} \|A v \|_H^{1/2} \| w \|_H, \quad u \in V, v \in D(A), w \in H, \\[2pt]
    |b(u, v, w)| &\leq C_G \|u\|_H \| v \|_V
    \| w \|_H^{1/2} \|A w \|_H^{1/2}, \quad \quad \quad \quad \, u \in H, v \in V, w \in
    D(A).
    \end{split}
\eeq Moreover, \beq \bl{Buv}
    \begin{split}
    \inptv{\Pi(u, v), \ w} &= - \inptv{\Pi(u, w),\  v}, \quad \tup{for} \; u, v, w \in V, \\[2pt]
    \inptv{\Pi(u, v), \ v} &= 0, \quad \tup{for} \; u, v \in V.
    \end{split}
\eeq For $u \in D(A)$, we have $B(u) \in H$,
\beq \bl{Bu}
    \hn{B(u)} \leq C_G \hn{u}^{1/2} \vn{u} \hn{Au}^{1/2},
\eeq
and especially, for the periodic case, we have
\begin{equation*}
    \inptv{\Pi(v, v),\  \Delta v}=0,\quad \tup{for} \;  v\in D(A) \textrm{  (c.f. \cite[Lemma 3.1, Page 19]{rT95})}.
\end{equation*}
\end{lemma}

The proof of Lemma \ref{L:bprop} can be found in \cite{rT84, rT88}. The first inequality in \eqref{buvw} will be most useful in this work and
its coefficient equals $2^{1/2}$ for any bounded and locally smooth domain in space dimension $n = 2$, which was proved in \cite[Lemma
3.4]{rT84}. The following lemma shows the regularity of functions in $H_0^1 (G)$ for a 2D domain $G$, whose proof is available in~\cite{rT84}.

\begin{lemma}
    For any two-dimensional open set $G$, we have
\beq \bl{2dreg}
    \| v \|_{L^4 (G)} \leq 2^{1/4} \| v \|_{L^2 (G)}^{1/2} \, \| \nabla v \|_{L^2 (G)}^{1/2},\ v \in H_0^1 (G).
\eeq

\end{lemma}

We have the following two versions of Gronwall-Bellman inequalites, whose proof is referred to \cite{DE92,FR75}.

(The Gronwall-Bellman Inequality): If a nonnegative scalar function $g(t)$
is continuous on $[0, T]$ and satisfies
\begin{equation} \label{Gwl}
    g (t) \leq (\geq)g(T) +  \int_t^T (\alpha g(s)+h(s))\, ds, \quad t \in [0,
    T],
\end{equation}
where $ \alpha \geq 0$ is a  constant and $h: [0,T]\rightarrow \mathbb{R}$
is integrable, then
\begin{equation} \label{Gwln}
    g(t) \leq(\geq)  \, e^{\alpha (T - t)}g(T)+\int_t^Te^{\alpha (s-t)}h(s)\  ds, \quad t \in [0, T].
\end{equation}

(The Stochastic Gronwall-Bellman Inequality): Let $(\Omega,\mathcal{F},\mathbb{F},P)$ be a filtered probability space whose filtration
$\mathbb{F}=\{ \mathcal{F}_t:t\in [0,T] \}$ satisfies the usual conditions. Suppose $\{ Y_s\}$ and $\{ X_s\}$ are optional integrable
processes and $\alpha$ is a nonnegative constant. If for all $t$, the map $s\mapsto E[Y_s|\mathcal{F}_t]$ is continuous almost surely and
$$Y_t\leq(\geq) E\left[\int_t^T(X_s+\alpha Y_s)\,ds+Y_T\,\Bigm |\,\mathcal{F}_t\right],$$ then we have almost surely
$$ Y_t\leq (\geq)e^{\alpha(T-t)}E[Y_T|\mathcal{F}_t]+E\left[ \int_t^Te^{\alpha (s-t)}
 X_s \,ds\, \Bigm |\,\mathcal{F}_t \right],\quad \forall t\in[0,T].   $$

In this paper, we prove the existence and uniqueness of an adapted solution to the terminal value problem \eqref{oreq} of a two-dimensional
backward stochastic Navier-Stokes equation with nonlinear forcing and the random perturbation driven by the Brownian motion.   We use the
spectral approximation, combined with the truncation and variational techniques, which is also a kind of compactness method.
The methodology features an interactive analysis based
on the regularity of the deterministic Navier-Stokes dynamics and the
stochastic properties of the It\^{o}-type diffusion processes.


\section{\textbf{Spectral approximations and estimates}}


In this section we consider the spectral approximation of the problem
\eqref{pb} obtained by orthogonally projecting the equation and the
terminal data on the finite dimensional space
$$
    H_N = \tup{Span} \, \{e_1, e_2, \cdots, e_N \}.
$$
 Define $$P_N: V' \to H_N, P_N f = \sum_{i=1}^{N} \lag f,\ e_i \rag
e_i,\quad f\in V'.$$ Then $\|P_Nf\|^2=\sum_{i=1}^N |\lag f,e_i \rag|^2$
 and $P_N$ is the orthogonal projection on $H_N$, which is called the spectral
projection. It is worth noting that $\|\cdot\|$ and $\|\cdot\|_V$ are
equivalent in $H_N$ and that $H_N=V_N:=P_N V$. Define \beq
    \begin{split}
    A^N = P_N A, \; B^N (u) = P_N B (u),\;
    J^NZ=P_N{JZ}:=\sum_{i=1}^{N}\lag JZ, e_i \rag e_i;  \; \\
    f^N (\cdots) = P_N f (\cdots), \; Z^N (t) = P_N Z(t), \; \tup{and} \; \xi^N = P_N \xi.
\end{split}
\eeq
 Then the projected, $N$-dimensional problem of approximation to the
problem \eqref{pb} is defined to be
\begin{equation}\label{appb}
  \left\{\begin{array}{l}
    \begin{split}
     d u^N (t) =&
        \left(\nu A^N u^N (t) - B^N (u^N (t))-J^NZ^N(t)
        -  f^N (t,u^N (t), Z^N (t))\right)\, dt\\
    &
        + Z^N (t) \, dW_t, \quad t \in [0,T); \\
        u^N(T) =
        & \xi^N.
    \end{split}
  \end{array}\right.
\end{equation}

 Note that the projection does not affect the Brownian motion
$\{W_t\}_{t \geq 0}$, and also that, the finite dimensional approximation
equation \eqref{appb} does not satisfy the conditions listed in
\cite{Hu-Briand-Pardoux03}.

We shall conduct \emph{a priori} estimates for the  adapted solution to the
finite dimensional approximation problem \eqref{appb}.

First, by means of Young's inequality
$$ab\leq \frac{1}{p}a^p \varepsilon ^p
+\frac{1}{q\varepsilon^q}b^q,\quad \frac{1}{p}+\frac{1}{q}=1, \quad ab\geq 0,~\varepsilon>0,  $$ under the Assumptions (A1)-(A3), we have
\begin{equation}\label{coercivity inequality}
    \begin{split}
    &2\inptv{\Phi(t,\phi,\varphi),\ \phi}- \| \varphi \|^2\\
    =&\ 2\inptv{-\nu A\phi+B(\phi)+J\varphi+f(t,\phi,\varphi),\ \phi}- \| \varphi \|^2\\
    =&- 2\nu \inptv{A \phi
      ,\  \phi}  - 2 \lag f (t, \phi, \varphi),\  \phi\rag
       -2\lag \varphi,\ (\sigma\cdot\nabla) \phi \rag
      - \| \varphi \|^2 \\
    \leq&
        -2\nu \|\phi\|_V^2+2 ( g(t)+\epsilon \|\varphi\|^2+\varrho(\epsilon) \| \phi\|^2+
      \beta\|\phi\|_V\|\phi\|)\\
     & +2\frac{1}{\bar{\lambda}} \|\varphi\|\|\bar{\lambda}(\sigma\cdot \nb) \phi\|
       - \| \varphi \|^2\quad
    \textrm{(choose $\epsilon$ small enough)}\\
    \leq&
    -2\lambda\|\phi\|_V^2-\frac{\bar{\lambda}^2-1}{2\bar{\lambda}^2}\|\varphi\|^2+
    2g(t)+\frac{\bar{\lambda}^2-1}{4\bar{\lambda}^2}\|\varphi\|^2
    +\lambda \|\phi\|_V^2+C\|\phi\|^2\\
    =&
    -\lambda\|\phi\|_V^2-\frac{\bar{\lambda}^2-1}{4\bar{\lambda}^2}\|\varphi\|^2+
    2g(t)+C\|\phi\|^2,
     \ (\phi,\varphi)\in V\times H,
    \end{split}
\end{equation}
where the constant $C$ depends only on $\lambda,\bar{\lambda},\varrho$ and
$\beta$.

\begin{lemma} \bl{L:apest}
    Let the conditions of Theorem \ref{thm main} hold. If $(u^N (\cdot), Z^N (\cdot))\in \ms{M}$ is an adapted solution of the
    problem \eqref{appb}, then
we have almost surely

(1). \begin{equation}\label{apriori 01}
  \begin{split}
    &\sup_{t\in [0,T]}\left\{\|u^N(t)\|^2
     +E_{\mathcal{F}_t}\left[ \int_t^T \|u^N(s)\|_V^2+ \|Z^N(s)\|^2\  ds \right]\right\}\\
     \leq&\  C\left( \|g\|_{L^\infty(\Omega, L^1([0,T]))}
         +\|\xi\|^2_{L^\infty(\Omega,H)}  \right),
  \end{split}
\end{equation}
where $C$ is a constant depending only on
$T,\nu,\lambda,\bar{\lambda},\beta$ and $\varrho$;

(2). \begin{equation}\label{apriori 02}
  \begin{split}
    &\sup_{t\in [0,T]}\left\{\|u^N(t)\|_V^2
     +E_{\mathcal{F}_t}\!\left[ \int_t^T\!\! \|Au^N(s)\|^2\!+\!\|Z^N(s)\|_V^2 ds \right]\right\}
    \\
        \leq&\
         C\left(\|g\|_{L^\infty(\Omega,L^1([0,T]))}
         +\|\xi\|^2_{L^\infty(\Omega,V)} \right)
  \end{split}
\end{equation}
with $C$ being a constant depending only on
$\nu,\lambda,\bar{\lambda},\beta,\varrho,\rho_1$ and $T$.
\end{lemma}

\bpr Applying the backward It\^{o} formula to the scalar-valued, stochastic process $\| u^N
    (t) \|^2$, and noting that
$$
    \inptv{B^N (u^N (t)), u^N (t)} = 0,
$$ we have
\begin{equation*}
    \begin{split}
    &\| u^N (t)\|^2\\
    =&\  \|\xi^N \|^2 - 2\nu \int_t^T\!\!\!\! \inptv{A^N u^N
      (s),\, u^N (s)} \,ds + 2 \int_t^T\!\!\!\! \lag f^N (s, u^N (s), Z^N
      (s)),\, u^N (s)\rag \,ds \\
      & +2\int_t^T\!\!\!\! \inptv{J^NZ^N(s),\,u^N(s)}\,ds- 2\int_t^T\!\!\!\! \lag Z^N (s),\, u^N (s) \rag \, dW_s
      - \int_t^T\!\!\!\!\| Z^N (s) \|^2 \, ds\\
    =&\ \|\xi^N \|^2 - 2\nu \int_t^T\!\!\!\! \inptv{A u^N
      (s),\, u^N (s)} \,ds + 2 \int_t^T\!\!\!\! \lag f (s, u^N (s), Z^N
      (s)),\, u^N (s)\rag \,ds \\
      & +2\int_t^T\!\!\!\! \inptv{JZ^N(s),\, u^N(s)}\,ds- 2\int_t^T\!\!\!\! \lag Z^N (s),\, u^N (s) \rag \, dW_s
      - \int_t^T\!\!\!\!\| Z^N (s) \|^2 \, ds.
    \end{split}
\end{equation*}
In view of~\eqref{coercivity inequality}, we have
   \begin{equation}\label{Ito formula for square}
    \begin{split}
    &\| u^N (t)\|^2\\
    =&\ \|\xi^N\|^2+\int_t^T \left(2\lag\Phi(s,u^N(s),Z^N(s)),\ u^N(s)    \rag -\|Z^N(s)\|^2 \right)\,ds\\
    &\ - 2\int_t^T \lag Z^N (s), \ u^N (s)  \rag \, dW_s\\
    \leq \ &
        \|\xi^N\|^2- 2\int_t^T \lag Z^N (s),\  u^N (s)  \rag \, dW_s \\
        &+\int_t^T\left(-\lambda \|u^N(s)\|_V^2
        -\frac{\bar{\lambda}^2-1}{4\bar{\lambda}^{2}}\|Z^N(s)\|^2
        +2g(s)+C\|u^N(s)\|^2   \right)\,ds
    \end{split}
\end{equation}
where the constant $C$ is independent of $N$.
 Since
 \begin{equation}\label{lemma2.1 2}
   \begin{split}
    &E\left[\sup_{\tau\in[t,T]}\left|\int_{\tau}^T \lag Z^N (s),\  u^N (s) \rag
        \ dW_s\right|\right]\\
   \leq &\
    2E\left[\sup_{\tau\in [t,T]} \left|\int_t^{\tau} \lag Z^N (s),\  u^N (s) \rag \ dW_s\right| \right] \\
   \leq &\  C
    E\left[ \left(\int_t^T \|Z^N(s)\|^2\|u^N(s)\|^2\ ds \right)^{1/2}\right]\quad
    (\textrm{by BDG inequality}) \\
    \leq &\
     CE\left[\sup_{s\in [t,T]}\|u^N(s)\|\left( \int_t^T
     \|Z^N(s)\|^2\ ds\right)^{1/2}\right]\\
     \leq &\
    (1/2)E\left[\sup_{s\in[t,T]}\|u^N(s)\|^2\right]
        +CE\left[ \int_t^T
     \|Z^N(s)\|^2\ ds \right],
    \end{split}
 \end{equation}
taking the conditional expectation on both sides of the second inequality
of \eqref{Ito formula for square}, we obtain
\begin{equation*}
  \begin{split}
    &\|u^N (t)\|^2
    +\lambda \ef{t}\left[\int_t^T\|u^N(s)\|^2_V\ ds\right]
    +\frac{\bar{\lambda}^2-1}{4\bar{\lambda}^2}\ef{t}\left[    \int_t^T\|Z^N(s)\|^2\ ds\right] \\
   \leq &\
    \ef{t}\left[ \| \xi \|^2 \right]+ C\ef{t}\left[
     \int_t^T\left(g(s)+\|u^N(s)\|^2\right)\ ds \right],\ a.s..
 \end{split}
\end{equation*}
From the stochastic Gronwall-Bellman inequality, it follows that
\begin{equation}\label{apriori 1}
  \begin{split}
    &\sup_{t\in [0,T]}\left\{\|u^N(t)\|^2
     +E_{\mathcal{F}_t}\left[ \int_t^T \|u^N(s)\|_V^2\ ds \right]
        +E_{\mathcal{F}_t}\left[\int_t^T \|Z^N(s)\|^2 \ ds \right]\right\}\\
        &\leq C\left(  \|g\|_{L^\infty(\Omega, L^1([0,T]))}
         +\|\xi^N\|^2_{L^\infty(\Omega,H)} \right)\\
        &\leq C\left( \|g\|_{L^\infty(\Omega, L^1([0,T]))}
         +\|\xi\|^2_{L^\infty(\Omega,H)}  \right),\ a.s.
  \end{split}
\end{equation}
where $C$ is a constant depending only on
$T,\nu,\lambda,\bar{\lambda},\beta$ and $\varrho$.

On the other hand, as $(B(u),\ \Delta u)=0$, using It\^o formula, we have
\begin{equation*}
  \begin{split}
    \|u^N(t)\|_V^2
      =&\  \|\xi^N \|_V^2 - 2\nu \int_t^T \inptv{A^N u^N
      (s),\  A^Nu^N (s)}\, ds \\
      &+2 \int_t^T \lag f^N (s, u^N (s), Z^N
      (s)),\  A^Nu^N (s)\rag\, ds  \\
        &+2\int_t^T \inptv{J^NZ^N(s),\ A^Nu^N(s)}\, ds\\
      &- 2\int_t^T \lag Z^N (s),\  A^Nu^N (s) \rag \, dW_s
      - \int_t^T\| Z^N (s) \|_V^2 \, ds\\
    =&\ \|\xi^N \|_V^2 - 2\nu \int_t^T \!\!\!\|A u^N(s)\|^2\, ds + 2 \int_t^T \!\!\!\lag f (s, u^N (s), Z^N
      (s)),\  Au^N (s)\rag\, ds\\
      & - \int_t^T\| Z^N (s) \|_V^2 \, ds  -2\int_t^T\sum_{i=1}^2 \lag  \nabla (Z^N)^i(s),\  (\sigma\cdot \nabla) \nabla (u^N)^i(s)) \rag \, ds\\
      & - 2\int_t^T \lag Z^N (s),\  Au^N (s) \rag \, dW_s,\ t\in[0,T].
  \end{split}
\end{equation*}
By the integration-by-parts formula and the fact that the integrals on the boundary $\partial G$ of $G$ vanish, we obtain
\begin{equation*}
  \begin{split}
    \|A\phi\|^2=\sum_{i=1}^2\|\nabla \phi^i\|_V^2,\ \ \forall\,\phi\in D(A).
  \end{split}
\end{equation*}
It follows that
\begin{equation*}
  \begin{split}
    &
        - 2\nu \int_t^T \!\!\!\|A u^N(s)\|^2\, ds
        -2\int_t^T\sum_{i=1}^2 \lag  \nabla (Z^N)^i(s),\  (\sigma\cdot \nabla) \nabla (u^N)^i(s)) \rag \, ds\\
    \leq&\
        -2\nu\sum_{i=1}^2\int_t^T\!\!\!\|\nabla (u^N(s))^i\|_V^2\,ds
        +\int_t^T \sum_{i=1}^2\|\bar{\lambda}(\sigma\cdot\nabla)\nabla (u^N)^i(s)\|^2\,ds
        \\
    &\
        +\bar{\lambda}^{-2}\int_t^T\!\!\!\|Z^N(s)\|_V^2\,ds\\
    \leq&\
        -2\lambda\sum_{i=1}^2\int_t^T\!\!\!\|\nabla (u^N(s))^i\|_V^2\,ds
        +\bar{\lambda}^{-2}\int_t^T\!\!\!\|Z^N(s)\|_V^2\,ds,\ a.s., \forall \,t\in[0,T].
  \end{split}
\end{equation*}
Therefore, we have
\begin{equation}\label{ito formula V}
    \begin{split}
    &\|u^N(t)\|_V^2\\
    \leq&\
        \|\xi \|_V^2- 2\int_t^T \lag Z^N (s),\  Au^N (s) \rag \, dW_s
        -\frac{\bar{\lambda}^2-1}{2\bar{\lambda}^2}\int_t^T\| Z^N (s) \|_V^2 \, ds\\
        & - 2\lambda \int_t^T \|A u^N(s)\|^2 \,ds+ 2 \int_t^T \lag f (s, u^N (s), Z^N(s)),\  Au^N (s)\rag \,ds \\
    \leq&\
        \|\xi \|_V^2- 2\int_t^T \lag Z^N (s),\  Au^N (s) \rag \, dW_s
        -\frac{\bar{\lambda}^2-1}{2\bar{\lambda}^2}\int_t^T\| Z^N (s) \|_V^2 \, ds\\
        &\ - \lambda \int_t^T \|A u^N(s)\|^2 ds + \frac{1}{\lambda} \int_t^T \|f (s, u^N (s), Z^N(s))\|^2\, ds \\
    \leq&\
        \|\xi \|_V^2- 2\int_t^T \lag Z^N (s),\  Au^N (s) \rag \, dW_s
        -\frac{\bar{\lambda}^2-1}{2\bar{\lambda}^2}\int_t^T\| Z^N (s) \|_V^2 \, ds\\
        &\
        - \lambda \int_t^T\!\!\! \|A u^N(s)\|^2 ds+ \frac{1}{\lambda} \int_t^T \left[g(s)
            +\beta(\|u^N(s)\|_V^2+\|Z^N(s)\|^2)\right]\rho_1(u^N)\, ds \\
    \leq&\
        \|\xi\|_V^2- 2\int_t^T \lag Z^N (s),\  Au^N (s) \rag \, dW_s
        -\frac{\bar{\lambda}^2-1}{2\bar{\lambda}^2}\int_t^T\| Z^N (s) \|_V^2 \, ds\\
        &\ - \lambda \int_t^T \|A u^N(s)\|^2 \,ds
        + C \int_t^T \left[g(s)
            +\beta(\|u^N(s)\|_V^2+\|Z^N(s)\|^2)\right] \,ds.
  \end{split}
\end{equation}
As
 \begin{equation}
   \begin{split}
    &E\left[\sup_{\tau\in[t,T]}\left|\int_{\tau}^T \lag Z^N (s), Au^N (s) \rag
        dW_s\right|\right]\\
   \leq &\
    2E\left[\sup_{\tau\in [t,T]} \left|\sum_{i=1}^2\int_t^{\tau} \lag \nabla (Z^N)^i (s),\nabla (u^N)^i (s) \rag dW_s\right| \right] \\
   \leq &\  C
    E\left[ \left(\int_t^T \|Z^N(s)\|_V^2\|u^N(s)\|_V^2ds \right)^{1/2}\right]\quad \quad
    (\textrm{by BDG inequality}) \\
     \leq &\
    (1/2)E\left[\sup_{s\in[t,T]}\|u^N(s)\|_V^2\right]
        +CE\left[ \int_t^T
     \|Z^N(s)\|_V^2ds \right]\\
     \leq &\
     C(N)\left\{ E\left[\sup_{s\in[t,T]}\|u^N(s)\|_V^2\right]
        +E\left[ \int_t^T
     \|Z^N(s)\|_V^2ds \right]  \right\},
    \end{split}
 \end{equation}
taking conditional expectation on both sides of \eqref{ito formula V}, we
get
\begin{equation}
  \begin{split}
    &\|u^N (t)\|_V^2
    +\lambda \ef{t}\left[\int_t^T\|Au^N(s)\|^2\,ds\right]
    +\frac{\bar{\lambda}^2-1}{2\bar{\lambda}^2}\ef{t}\left[
                            \int_t^T\|Z^N(s)\|_V^2\,ds\right] \\
   \leq &\
    \ef{t}\left[ \| \xi \|_V^2 \right]+ C\ef{t}\left[
     \int_t^T(g(s)+\|u^N(s)\|_V^2+\|Z^N(s)\|^2)\,ds \right].
  \end{split}
\end{equation}
In view of \eqref{apriori 1}, we conclude that, with probability 1,
\begin{equation}\label{apriori 2}
  \begin{split}
    &\sup_{t\in [0,T]}\left\{\|u^N(t)\|_V^2
     +E_{\mathcal{F}_t}\left[ \int_t^T \|Au^N(s)\|^2ds \right]
        +E_{\mathcal{F}_t}\left[\int_t^T \|Z^N(s)\|_V^2 ds \right]\right\}\\
        &\leq C\left( \|g\|_{L^\infty(\Omega,L^1([0,T]))}
         +\|\xi\|^2_{L^\infty(\Omega,V)} \right),
  \end{split}
\end{equation}
where $C$ is a constant depending only on
$\nu,\lambda,\bar{\lambda},\beta,\varrho,\rho_1$ and $T$. \epr

\begin{lemma} \bl{L:difest}
For any $u, v \in V$ and $\phi,\varphi\in H$,

 \beq \bl{umv}
    |\inptv{B(u) - B(v), u - v}| \leq \frac{\lambda}{4} \|u - v\|_V^2
    + \frac{2}{\lambda} \|v \|_V^2 \| u - v \|^2.
\eeq Moreover, under Assumptions (A1)-(A3),  there exists a positive constant $K$ depending on $\bar{\lambda}$ such that
 \beq \bl{coep}\begin{split}
     &-2\inptv{\Phi(t,u,\phi)-\Phi(t,v,\varphi), w}\\
    + & \|w \|^2\left(K+\frac{4}{\lambda}\|v\|_V^2+K\rho^2(v)\right)+
    \frac{\bar{\lambda}^2+1}{2\bar{\lambda}^2}\|\bar{w}\|^2
      \geq \frac{\lambda}{2} \| w \|_V^2,~\quad
  \end{split}
\eeq
 holds almost surely for any $t\in [0,T],$ $u, v \in V$ and $\varphi,\phi\in H$ with $w := u - v, \bar{w}:=\phi-\varphi$,
  and $\Phi$ being defined by~\eqref{lemma2.2 0}. Define
$$
    r_1(t) =  \int_0^t \left(K+\frac{4}{\lambda}\| u(s)\|_V^2+K\rho^2(u(s))\right) \, ds
    $$
and $$ r_2 (t) =  \int_0^t\left(K+\frac{4}{\lambda}\| v(s)\|_V^2
    +K\rho^2(v(s))\right)
    \, ds,
$$
for arbitrary $u, v \in L_{\ms{F}}^2 (\gw; L^2 (0, T; V))$, and let $w(\cdot) = u(\cdot) - v(\cdot)$. Then for any $\phi,\varphi \in
L_{\ms{F}}^2 (\gw; L^2 (0, T; H)) $ and $\bar{w}(\cdot):=\phi(\cdot)-\varphi(\cdot),$ we have
 for $i=1,2$
 \beq \bl{req}
    -\inptv{2\Phi(t,u,\phi)-2\Phi(t,v,\varphi) + \frac{d r_i(t)}{dt}
    w, w}+   \frac{1+\bar{\lambda}^2}{2\bar{\lambda}^2}
        \|\bar{w}\|^2 \geq 0, \ a.s..
\eeq
\end{lemma}

\bpr Let $w = u - v$. Then
\begin{align*}
    &\inptv{B(u) - B(v), u - v}\\
    = &-\inptv{\Pi (u,w), u} + \inptv{\Pi (v,w), v} \\[2pt]
    = &- \inptv{\Pi (u, w), v} + \inptv{\Pi (v, w), v} = - \inptv{B(w), v}.
\end{align*}
By the first inequality in Lemma \ref{L:bprop} we can get
\begin{align*}
    |\lag & B(u) - B(v), u - v \rag_{V',V}|  = |\lag B(w), v \ragv | = |\lag \Pi (w,v), w \rag_{V',V}| \\[2pt]
    & \leq 2^{1/2} \|u - v\| \|u - v\|_V \| v \|_V \leq \frac{\lambda}{4} \|u - v\|_V^2
    + \frac{2}{\lambda} \|u - v\|^2
    \|v\|_V^2.
\end{align*}
It follows from the Assumptions (A1)-(A3) that
\begin{equation}\label{remark using}
  \begin{split}
    &-2\nu\inptv{Aw,w}+2\inptv{J\bar{w},w}+2\lag
        f(t,u,\phi)-f(t,v,\varphi),w \rag  \\
    \leq&
        -2\nu\|w\|_V^2+2\lag\bar{w},\ (\sigma\cdot\nb) w
        \rag+2\rho(v)\|w\|^2+2\rho(v)\|w\|(\|w\|_V+\|\bar{w}\|) \\
    \leq&
        -\frac{3\lambda}{2}\|w\|_V^2+\frac{1}{\bar{\lambda}^2}\|\bar{w}\|^2
        +(1-\frac{1+\bar{\lambda}^2}{2\bar{\lambda}^2})\|\bar{w}\|^2\\
        &+(C(\bar{\lambda})\rho^2(v)+2\rho(v))\|w\|^2\\
    \leq&
        -\frac{3\lambda}{2}\|w\|_V^2+
            \frac{1+\bar{\lambda}^2}{2\bar{\lambda}^2}
                \|\bar{w}\|^2
                    +(K+K\rho^2(v))\|w\|^2,
  \end{split}
\end{equation}
where the constant $K$ only depends on $\bar{\lambda}$. Hence, in view of \eqref{umv}, we obtain \eqref{coep}.

Then \eqref{req} for $i = 2$ follows from \eqref{coep} by direct calculation. The case of $i = 1$ in \eqref{req} is shown in a similar way.
\epr

\section{\textbf{Solutions of the finite dimensional systems}}
In this section, we consider the existence of an adapted solution to the
projected, $N$-dimensional problem \eqref{appb} of the 2D backward
stochastic Navier-Stokes equations which we also call the
finite dimensional system. To solve the finite dimensional system
\eqref{appb}, we shall  make use of the result of Briand et al.
\cite{Hu-Briand-Pardoux03}.

Consider the following backward stochastic differential equation (BSDE in short):
\begin{equation}\label{BSDE}
  \begin{split}
    Y(t)=\zeta+\int_t^Tg(s,Y(s),q(s))\, ds-\int_t^Tq(s)\, dW_s,
  \end{split}
\end{equation}
where $\zeta$ is an $\mathbb{R}^N$-valued $\mathscr{F}_T$-measuable random vector and the random function
$$g:~[0,T]\times \Omega\times
\mathbb{R}^N\times\mathbb{R}^N\rightarrow \mathbb{R}^N$$ is $\mathcal {P}\times \mathscr{B}(\mathbb{R}^N)\times
\mathscr{B}(\mathbb{R}^N)$-measurable.

The following lemma comes from \cite[Theorem 4.2]{Hu-Briand-Pardoux03}.
\begin{lemma}\label{BSDE lem}
  Assume that $g$ and $\zeta$ satisfy the following four conditions.

  (C1). For some $p>1$, we have
  $$ E\left[|\zeta|^p+\left( \int_0^T|g(t,0,0)|\ ds\right)^p\right]<\infty. $$

  (C2).
   There exist constants $\alpha \geq 0$ and $\mu\in\mathbb{R}$ such
  that almost surely we have  for each $(t,y,y',z,z')\in[0,T]\times\mathbb{R}^N\times\mathbb{R}^N\times\mathbb{R}^N
  \times \mathbb{R}^N,$
  \begin{equation}\label{BSDE lem C2 1}
    |g((t,y,z)-g(t,y,z')|\leq \alpha|z-z'|,
  \end{equation}
  \begin{equation}\label{BSDE lem C2 2}
    \lag y-y',g(t,y,z)-g(t,y',z)\rag \leq \mu |y-y'|^2\quad \hbox{ \rm (monotonicity condition).}
  \end{equation}

  (C3).  The function $y\mapsto
  g(t,y,z)$ is continuous for any $(t,z)\in [0,T]\times\mathbb{R}^N$.

  (C4). For any $r>0$, the random process\\ $$\left\{\psi_r(t):=\sup_{|y|\leq r}|g(t,y,0)-g(t,0,0)|,\  t\in [0,T]\right\}$$
   lies in the space $L^1(\Omega\times [0,T])$.
Then BSDE~\eqref{BSDE} admits a unique solution  $(Y,q)\in L^p(\Omega,C([0,T],\mathbb{R}^N))\times \L^p(\Omega,L^2([0,T],\mathbb{R}^N)).$
\end{lemma}
\begin{remark}
  It is worth noting that our finite dimensional system does not satisfy the
 monotonicity condition (C2). In fact, by Lemma \ref{L:difest} our
 finite dimensional system only satisfies a local monotonicity condition in
 some sense, which prevents us to directly use this lemma to our finite dimensional
 system.
\end{remark}

\begin{lemma} \bl{L:trn}
    For any $M,N\in \mathbb{Z}^+$, define the function of
truncation $R_M (\cdot)$ to be a $C^2$ function on $H_N$ such that for
$X=\sum_{i=1}^N x_i e_i,$
$$
    R_M (X) =
    \begin{cases}
      1,  &\textup{if} \; \|X\| \leq M; \\[2pt]
      \tup{$\in (0,1)$,} &\textup{\rm if} \; M <
      \|X\| < M + 1; \\[2pt]
      0, &\textup{if} \; \|X\| \geq M + 1.
    \end{cases}
$$
Thus $R_M (\cdot)$ is uniformly Lipschitz continuous. For each
$n\in\mathbb{Z}^+$, denote $\varphi_n(z)=zn/(\|z\|\vee n)$, $z\in H^N$  and
set
$$
\Phi^{N,M,n}(t,y,z)=R_M(y)\frac{n}{h_{M}(t)\vee n}
  P_N\Phi(t,y,\varphi_n(z)),
$$
where
\begin{equation}
    \begin{split}
         h_{M}(t)&=4\left\{\left(g(t)+\beta C_N (M+1)^2\right)\esssup_{\|w\|\leq M+1}|\rho_1(w)|
                    +C_{N,\nu}(M+1)^2 \right\}^{1/2}\\
                 &\geq
                 \esssup_{\|w\|\leq M+1}\|\Phi(t,w,0)\|
    \end{split}
\end{equation}
and $h_M\in L^1(\Omega\times [0,T]).$
 Then under
Assumptions (A1)-(A3), $\Phi^{N,M,n}$ satisfies the conditions (C2)-(C4) in Lemma \ref{BSDE lem}.
\end{lemma}

\begin{proof}
Under Assumptions (A1)-(A3) and Remark \ref{rmk Lip wrt Z}, we only need verify \eqref{BSDE lem C2 2}, i.e.,
 there is a uniform constant $C_{N,M,n} > 0$ such that
\beq \label{finite-dim monotone}
    \lag \Phi^{N,M,n} (t,X,Z) - \Phi^{N,M,n} (t,Y,Z),\ X-Y\rag \leq
    C_{N,M,n} \| X - Y \|^2,\quad a.s.,
\eeq for any $X, Y, Z \in H_N$ and all $t\in [0,T]$.
 For any $X, Y \in H_N$,  inequality~\eqref{finite-dim monotone} holds trivially if $\|X\|>M+1 $ and $\|Y\|>M+1$.
 Thus, it is sufficient to consider the case of $\|Y\|\leq M+1$.  We  have
\begin{equation}
  \begin{split}
    &\lag \Phi^{N,M,n} (t,X,Z) - \Phi^{N,M,n} (t,Y,Z),\ X-Y\rag\\
    =&\
      R_M (X )\frac{n}{h_{M}(t)\vee n}
       \lag  \Phi (t,X,\varphi_n(Z)) - \Phi (t,Y,\varphi_n(Z)),\ X-Y\rag \\
     & +\frac{n}{h_{M}(t)\vee n}(R_M(X)-R_M(Y))\lag \Phi^{N,M,n}(t,Y,\varphi_n(Z)),\ X-Y \rag \\
   & \quad\quad \quad (~\textrm{by \eqref{coep} of Lemma \ref{L:difest}}) \\
    \leq &\
        \left(K+\frac{4}{\lambda}\|Y\|^2+K\rho^2(Y)\right)\|X-Y\|^2\\
        &\quad+C_M\|X-Y\|^2\frac{n}{h_{M}(t)\vee n}\|\Phi^{N,M,n}(t,Y,\varphi_n(Z))\|  \\
    \leq &\
        \left(K+\frac{4}{\lambda}\|Y\|^2+K\rho^2(Y)\right)\|X-Y\|^2\\
          & \quad+ C_{M}\|X-Y\|^2\frac{n}{h_{M}(t)\vee n} (h_{M}(t)+ C_{N,M}\cdot n )  \\
    \leq &\
        C_{M,N,n}\|X-Y\|^2,
  \end{split}
\end{equation}
which completes the proof.
\end{proof}

 \bthm \bl{T:prj} Let Assumptions \tup{(A1)-(A3)} hold. For
any $\xi \in L_{\ms{F}}^\infty (\gw; V)$, the projected problem \eqref{appb} admits a unique adapted solution $(u^N (\cdot), Z^N (\cdot))\in
\ms{M}$ for each given positive integer $N$, which satisfies
\begin{equation}\label{finit-dim esti}
\|(u^N,Z^N)\|_{\ms{M}} \leq C\left\{1+E\left[\|\xi\|^2\right] \right\},
\end{equation}
where $C$ is a constant independent of $N$.
 \ethm
\begin{proof}
\tbf{Step 1}. Let us verify the uniqueness part. Suppose $(u^N,Z^N)$ and $(v^N,Y^N)$ are two solutions of the projected problem \eqref{appb}.
Note that the a priori estimates in Lemma \ref{L:apest} holds for both $(u^N,Z^N)$ and $(v^N,Y^N)$. Denote by $(U^N,X^N)$ the pair of
processes $(u^N-v^N,Z^N-Y^N)$. Define
$$ r(t):=\int_0^t \left[K+\frac{4}{\lambda}\|v^N(s)\|_V^2+K\rho^2(v^N(s))\right]\, ds.$$ An application
of It\^o formula and Lemma \ref{L:difest} yields
\begin{equation}\label{Thm finite-dim unique}
  \begin{split}
    &e^{r(t)}\|U^N(t)\|^2\\
    =&\int_t^Te^{r(s)} \biggl[2 \left\lag
        \Phi(s,u^N(s),Z^N(s))-\Phi(s,v^N(s),Y^N(s)), \, U^N(s)
            \right\rag \\
     &\quad      -\|X^N\|^2- \|U^N(s)\|^2\left(K+\frac{4}{\lambda}\|v^N(s)\|_V^2+K\rho^2(v^N(s))\right)\biggr]\,ds\\
    & -2\int_t^Te^{r(s)}\lag U^N(s),X^N(s) \rag dW_s\\
    \leq &
     -\frac{\bar{\lambda}^2-1}{2\bar{\lambda}^2}
        \int_t^T e^{r(s)}\|X^N(s)\|^2ds-
        2\int_t^Te^{r(s)}\lag U^N(s),X^N(s) \rag\, dW_s.
  \end{split}
\end{equation}
Taking conditional expectations on both sides, we get
    $$e^{r(t)}\|U^N(t)\|^2
    +\frac{\bar{\lambda}^2-1}{2\bar{\lambda}^2}E_{\mathcal{F}_t}\left[
            \int_t^T e^{r(s)}\|X^N(s)\|^2ds
                        \right]\leq 0,\quad a.s.,
        \textrm{ for any }t\in[0,T],$$
 which implies the uniqueness.

\tbf{Step 2}. For any $N,M,n\in \mathbb{Z}^+$, following Lemma \ref{L:trn}, we can verify that the pair $(\xi^N,\Phi^{N,M,n})$ satisfies the
conditions (C1)-(C4) in Lemma \ref{BSDE lem}. Hence, by Lemma \ref{BSDE lem} there exists a unique solution $(u^{N,M,n},Z^{N,M,n})\in \ms{M}$
to the following BSDE:
\begin{equation}
    u^{N,M,n}(t)=\xi^N+\int_t^T\!\!
    \Phi^{N,M,n}(s,u^{N,M,n}(s),Z^{N,M,n}(s))\, ds-\int_t^T\!\!Z^{N,M,n}(s)\, dW_s.
\end{equation}
In a similar way to Lemma \ref{L:apest}, we deduce that there exists a
positive constant $K_1$ which is independent of $N,M$ and $n$ such that
\begin{equation}\label{Thm finite-dim ex 1}
    \sup_{t\in[0,T]}\|u^{N,M,n}(t)\|+E\left[ \int_0^T\!\! \|Z^{N,M,n}(s)\|^2\, ds
        \right]\leq K_1,\quad a.s..
\end{equation}
Therefore, taking $M>K_1$, we have $R_M(u^{N,M,n}(s))\equiv 1$ and $(u^{N,M,n},Z^{N,M,n})$ is independent of $M$. Thus, we write
$(u^{N,n},Z^{N,n})$ instead of $(u^{N,M,n},Z^{N,M,n})$ below. Moreover, there exists a positive constant $K_2$ independent of $n$ such that
\begin{equation}
  \begin{split}
    K+\frac{4}{\lambda}\|u^{N,n}(t)\|_V^2+K\rho^2(u^{N,n}(t))
     &\leq K_2,\\
        \|\Phi(t,u^{N,n}(t),\phi_1)-\Phi(t,u^{N,n}(t),\phi_2)\|
     &\leq K_2
        \|\phi_1-\phi_2\|,\ \ dP\otimes dt\textrm{-almost},
  \end{split}
\end{equation}
holds for all $\phi_1,\phi_2\in H$ and $N,n\in \mathbb{Z}^+$.

For $j\in\mathbb{Z}^+$, set $(U^{N},X^{N})=(u^{N,n+j}-u^{N,n},Z^{N,n+j}-Z^{N,n})$. Applying It\^o formula similar to \eqref{Thm finite-dim
unique}, we get
\begin{equation}\label{thm finite-dim ex 2}
  \begin{split}
    &e^{K_2t}\|U^N(t)\|^2
     +
        \frac{\bar{\lambda}^2-1}{2\bar{\lambda}^2}\int_t^T
        e^{K_2s}\|X^N(s)\|^2\,ds\\
     \leq&\
     2\int_t^T\!\!\!\!e^{K_2s}\lag
        \Phi^{N,n+j}(s,u^{N,n}(s),Z^{N,n}(s))\!-\!\Phi^{N,n}(s,u^{N,n}(s),Z^{N,n}(s)),\ U^N(s)
            \rag ds
         \\
        &-2\int_t^T\!\!\!\!e^{K_2s}\lag U^N(s), X^N(s)   \rag\, dW_s\\
     (&\textrm{by \eqref{Thm finite-dim ex 1}})\\
     \leq&\
     4K_1 \int_t^T \!\!\!\! e^{K_2s}\|
        \Phi^{N,n+j}(s,u^{N,n}(s),Z^{N,n}(s))-\Phi^{N,n}(s,u^{N,n}(s),Z^{N,n}(s))
            \| \,ds\\
        &-2\int_t^T\!\!\!\!e^{K_2s}\lag U^N(s), X^N(s)   \rag\, dW_s.
  \end{split}
\end{equation}
On the other hand,
 \begin{equation}
   \begin{split}
    &E\left[\sup_{\tau\in[t,T]}|\int_{\tau}^T e^{K_2s} \lag U^N(s), X^N(s)
    \rag \,dW_s|\right]\\
   \leq & \ C
    E\left[ \left(\int_t^Te^{2K_2s} \|X^N(s)\|^2\|U^N(s)\|^2\,ds \right)^{1/2}\right]
    (\textrm{by BDG inequality}) \\
     \leq &\
    \epsilon E\left[\sup_{s\in[t,T]}(e^{K_2s}\|U^N(s)\|^2)\right]
        +C_{\epsilon}E\left[ \int_t^T
     \|X^N(s)\|^2e^{K_2s}\,ds \right],
    \end{split}
 \end{equation}
 with the positive constant $\epsilon$ to be determined later.
Then choosing $\epsilon$ to be small enough, we deduce from \eqref{thm finite-dim ex 2} that
\begin{equation*}
  \begin{split}
    &\|(U^N,X^N)\|_{\ms{M}}\\
    \leq&\  C E\left[ \int_0^T \!\!\!\|
        \Phi^{N,n+j}(s,u^{N,n}(s),Z^{N,n}(s))-\Phi^{N,n}(s,u^{N,n}(s),Z^{N,n}(s))
            \| \,ds  \right].
  \end{split}
\end{equation*}
As
\begin{equation}
  \begin{split}
    &\|
        \Phi^{N,n+j}(s,u^{N,n}(s),Z^{N,n}(s))-\Phi^{N,n}(s,u^{N,n}(s),Z^{N,n}(s))
            \|   \\
    \leq&\
       2K_2\|Z^{N,n}(s)\|\mathbb {I}_{\{\|Z^{N,n}(s)\|>n\}}\!+\!
        2K_2\|Z^{N,n}(s)\|\mathbb {I}_{\{h_{K_1}(s)>n\}}\!+\!
            2h_{K_1}(s)\mathbb {I}_{\{h_{K_1}(s)>n\}},
  \end{split}
\end{equation}
in view of \eqref{Thm finite-dim ex 1} and $h_{K_1}\in L^1(\Omega\times [0,T])$, we conclude that $(u^{N,n},Z^{N,n})$ is a Cauchy sequence in
$\ms{M}$. Denote the limit by $(u^N,Z^N)\in\ms{M}$. It is easily checked that $(u^N,Z^N)$ is a solution of the projected problem
\eqref{appb}.

\tbf{Step 3}. Estimate \eqref{finit-dim esti} follows from Lemma
\ref{L:apest}, which completes the proof.

\end{proof}

\section{\textbf{Proof of Theorem 2.1}}
\begin{proof}[Proof of Theorem \ref{thm main}]
  Our proof consists of the following four steps.

\tbf{Step 1}. By Theorem \ref{T:prj}, we have solved the projected problem
\eqref{appb} in $\ms{M}$. By Lemma \ref{L:apest}, we have
\begin{equation}\label{thm main prof step1 1}
  \begin{split}
    &
        \esssup_{(\omega,s)\in \Omega\times[0,T]}\|u^N(s)\|_V^2
        +E\left[ \int_0^T \|Au^N(s)\|^2+\|Z^N(s)\|_V^2 ds \right]
    \\
     \leq&\
      C\left( \|g\|_{L^\infty(\Omega,L^1([0,T]))}
         +\|\xi\|^2_{L^\infty(\Omega,V)} \right),
  \end{split}
\end{equation}
where $C$ is a constant depending only on
$\nu,\lambda,\bar{\lambda},\beta,\varrho,\rho_1$ and $T$
. Since we
get $ \|B(v)\|^2\leq C_G\|v\|\|v\|_V^2\|Av\|$ from Lemma \ref{L:bprop}, under Assumptions (A1)-(A3),
we conclude
\begin{equation*}
    E\left[ \int_0^T\left(\|B(u^N(s))\|^2
        +\|(\sigma\cdot\nabla) Z^N(s)\|^2+\|f(s,u^N(s),Z^N(s))\|^2 \right)\,ds \right]
        \leq C.
\end{equation*}
 Hence,
\begin{equation}
\|P_N\Phi(\cdot,u^N,Z^N)\|_{L^2(\Omega\times [0,T],H)}\leq
\|\Phi(\cdot,u^N,Z^N)\|_{L^2(\Omega\times [0,T],H)}\|\leq C.
\end{equation}
All the constants $C$s above are independent of $N$.

\tbf{Step 2}. Now we consider the weak convergence. Clearly,
$$
    \xi^N \to \xi \;\; \tup{strongly in} \; V,a.s., \textrm{ and }\|\xi^N\|_V \leq \|\xi\|_V \; \tup{as} \; N
    \to \infty,
$$
which implies that $\xi^N \rightarrow \xi $ in $L^p(\Omega,V)$ for any
$p\in (1,+\infty)$. Then the following weak and weak star convergence
results in respective spaces hold: there exists a subsequence
$\{N_k\}_{k=1}^\infty$ of $\{N\}$, such that, as $k \to \infty$, \beq
\bl{wkcv}
    \begin{split}
    u^{N_k} (\cdot) & \overset{w}{\longrightarrow} u(\cdot) \; \tup{in} \;
    \ms{L}_{\ms{F}}^2 (0, T; D(A))), \\
    u^{N_k} (\cdot) & \overset{w^*}{\longrightarrow} u(\cdot) \; \tup{in} \;
    L^\infty_{\ms{F}}
    (\gw; C ([0, T]; V)), \\
    Z^{N_k} (\cdot) & \overset{w}{\longrightarrow} Z(\cdot) \; \tup{in} \; \lf
    (\gw; L^2 (0, T; V)),\\
    \Phi(\cdot,u^{N_k},Z^{N_k}) (\cdot) & \overset{w}{\longrightarrow} \Gamma (\cdot) \; \tup{in} \; \lf
    (\gw; L^2 (0, T; H)), \\
    P_N\Phi(\cdot,u^{N_k},Z^{N_k}) (\cdot) & \overset{w}{\longrightarrow} \Psi(\cdot) \; \tup{in} \; \lf
    (\gw; L^2 (0, T; H)), \\
    \end{split}
\eeq where $u,Z,\Gamma$ and $\Psi$ are some functions in the respective
spaces.

 By the Burkh\"{o}lder-Davis-Gundy (BDG in short) inequality, we can get
\begin{equation}
    \begin{split}
    E &\left[\int_0^T \left\|\int_t^T Z^{N_k} (s)  dW_s \right\|_V^2 dt\right]
     \leq T
         E \left[\sup_{t \in [0, T]} \left\| \int_t^T Z^{N_k} (s) dW_s \right\|_V^2\right] \\
    &\leq 2T E \left[\sup_{t \in [0, T]} \left\| \int_0^t Z^{N_k} (s) dW_s \right\|_V^2\right]
    + 2 T E \left[ \left\| \int_0^T Z^{N_k} (s) dW_s\right \|_V^2\right] \\
    & \leq 4 T E \left[\sup_{t \in [0, T]} \left\| \int_0^t Z^{N_k} (s)dW_s \right\|_V^2\right]
     \leq 4L_1 T E \left[\int_0^T \|Z^{N_k} (s) |_V^2 ds \right]
 \end{split}
\end{equation}
 where $L_1 > 0$ is a uniform constant from the BDG inequality. Hence, as a bounded linear operator
on the space $\lf (\gw; L^2 (0, T; V))$, the mapping
$$
    \Upsilon: Z^{N_k} (\cdot) \longmapsto \int_{\cdot}^T Z^{N_k} (s) \,
    dW_s
$$
maps the weakly convergent sequence $\{Z^{N_k} (\cdot)\}$ to a weakly
convergent sequence $
    \left\{ \int_{\cdot}^T Z^{N_k} (s)  dW_s \right\}$
    $  \tup{in } \lf (\gw; L^2 (0, T; V))$ such that
\begin{equation}\bl{intzcv}
    \int_{\cdot}^T Z^{N_k} (s) \, dW_s \overset{w}{\longrightarrow} \int_{\cdot}^T Z(s) \,
    dW_s\textrm{ in } \lf (\gw; L^2 (0, T; V)), \quad \tup{as} \; \, k \to
    \infty.
\end{equation}
 Similarly it can be shown that, $\textrm{ as } k \to \infty$,
 \begin{equation} \bl{intf}
    \int_{\cdot}^T P_{N_k}\Phi (s, u^{N_k} (s), Z^{N_k} (s))\, ds
    \overset{w}{\longrightarrow} \int_{\cdot}^T \Psi(s)\, ds
     \textrm{ in } \lf (\gw; L^2 (0, T; H)).
\end{equation}
Define
\begin{equation}\label{thm main prof step31}
    \bar{u}(t)=\xi+\int_t^T \Psi(s) \ ds-\int_t^T Z(s) \ dW_s.
\end{equation}
 It is easily
checked that $\bar{u}=u$, $P\otimes dt$-almost. In view of \cite[Theorem 4.2.5]{PR07}, we conclude that $u\in
L^\infty(\Omega,C([0,T],V))$ and also,
    \begin{equation}\label{estimate}
    \begin{split}
       &
            \esssup_{(\omega,s)\in\Omega\times[0,T]}\|u(s)\|_V+ E\left[\int_0^T\!\!\! \|u(s)\|_{D(A)}^2 ds\!+\!\!
            \int_0^T \!\!\!\!\|Z(s)\|_V^2 ds  \right]
        \\
        \leq&\
            C\left\{\|g\|_{L^\infty(\Omega,L^1([0,T]))}+\|\xi\|_{L_{\ms{F}_T}^\infty (\gw; V)}^2 \right\},
    \end{split}
    \end{equation}
    where $C$ is a constant depending on
    $\nu,\lambda,\bar{\lambda},\beta,\varrho,\rho_1 \textrm{ and }
     T$.

 \tbf{Step 3}.
 For a notational convenience, we now use the index
  $N$ instead of $N_k$ for all the relevant subsequences.

As $\cup_{N=1}^\infty \ms{L}_{\ms{F}}^2 (0, T; P_N H)$ is dense in $\ms{L}_{\ms{F}}^2 (0, T; H)$ and it can be checked that $\Psi=\Gamma$ on
$\cup_{N=1}^\infty \ms{L}_{\ms{F}}^2 (0, T; P_N H)$, by a density argument we have $\Psi=\Gamma$. Thus, to show that $(u,Z)$ is a strong
solution of the 2D BSNSE problem \eqref{pb}, we only need prove
\begin{equation}\label{drift verification}
    \Psi(\cdot)=\Phi(\cdot,u,Z),\ a.s..
\end{equation}

For any $v\in L^\infty_{\ms{F}}(\Omega,C([0,T],V)),$ define
$$
r(t)=r(\omega,t):=\int_0^t(K+\frac{4}{\lambda}\|v(\omega,s)\|_V^2+K\rho^2(v(\omega,s)))\ ds,
\quad (\omega, t)\in \Omega\times[0,T],
$$
where the constant $K$ comes from \eqref{coep} in Lemma \ref{L:difest}. Applying It\^o formula
to compute $e^{r(t)}\|u^N(t)\|^2$, we have
\begin{equation*}
\begin{split}
  &E\left[ e^{r(t)}\|u^N(t)\|^2-e^{r(T)}\|u^N(T)\|^2 \right]  \\
  =&\
     E\bigg[\int_t^{T} e^{r(s)}\Big( \inptv{2P_N\Phi(s,u^N(s),Z^N(s)),u^N(s)}
     -\|Z^N(s)\|^2 \\
     &\,\,\,\,-\big(K+\frac{4}{\lambda}\|v(s)\|_V^2+K\rho^2(v(s))\big)\|u^N(s)\|^2\Big)\,ds\bigg]  \\
  =&\
    E\bigg[\int_t^{T} e^{r(s)}\bigg(
    2\inptv{\Phi(s,u^N(s),Z^N(s))-\Phi(s,v(s),Z(s)),u^N(s)-v(s)}\\
    &\quad -\|Z^N(s)-Z(s)\|^2
    -\left(K+\frac{4}{\lambda}\|v(s)\|_V^2+K\rho^2(v(s))\right)\|u^N(s)-v(s)\|^2
        \bigg)\,ds\bigg]\\
    &+E\bigg[\int_t^{T} e^{r(s)}\bigg(
        2\inptv{\Phi(s,u^N(s),Z^N(s))-\Phi(s,v(s),Z(s)),v(s)}\\
    &\quad+2\inptv{\Phi(s,v(s),Z(s)),u^N(s)}
    -2\lag Z^N(s), Z(s)\rag+\|Z(s)\|^2\\
    &\quad
        -\left(K+\frac{4}{\lambda}\|v(s)\|_V^2+K\rho^2(v(s))\right)\left( 2\lag u^N(s),v(s) \rag
        -\|v(s)\|^2\right)
       \bigg)\,ds\bigg]\\
  \leq&\
    E\bigg[\int_t^{T} e^{r(s)}\bigg(
        2\inptv{\Phi(s,u^N(s),Z^N(s))-\Phi(s,v(s),Z(s)),v(s)}\\
    &\quad+2\inptv{\Phi(s,v(s),Z(s)),u^N(s)} -2\lag Z^N(s), Z(s)\rag+\|Z(s)\|^2\\
    &\quad
        -\left(K+\frac{4}{\lambda}\|v(s)\|_V^2+K\rho^2(v(s))\right)\left( 2\lag u^N(s),v(s) \rag
        -\|v(s)\|^2\right)
       \bigg)\,ds\bigg].
\end{split}
\end{equation*}

Letting $N \rightarrow \infty,$ by Lemma \ref{L:difest} and the lower semicontinuity, we have for any nonnegative $\varphi\in
L^{\infty}(0,T),$
\begin{equation}\label{exsitence thm 1}
\begin{split}
  &E\left[\int_0^{T}\varphi(t)\left(e^{r(t)}
  \|u(t)\|^2-e^{r(T)}\|u(T)\|^2\right)dt
    \right]\\
  \leq&\
    \liminf_{N\rightarrow\infty}
        E\left[\int_0^T\varphi(t)\left(e^{r(t)}
        \|u^N(t)\|^2-e^{r(T)}\|u^N(T)\|^2\right)dt
            \right]\\
  \leq&\
    E\bigg[\int_0^{T}\varphi(t)\bigg(\int_t^{T} e^{r(s)}\bigg(
        2\inptv{\Psi(s)-\Phi(s,v(s),Z(s)),v(s)}\\
    &\quad+2\inptv{\Phi(s,v(s),Z(s)),u(s)}
     -2\lag Z(s), Z(s)\rag+\|Z(s)\|^2\\
    &\quad
        -\left(K+\frac{4}{\lambda}\|v(s)\|_V^2+K\rho^2(v(s))\right)\left( 2\lag u(s),v(s) \rag
        -\|v(s)\|^2\right)
       \bigg)ds\bigg)dt\bigg],
\end{split}
\end{equation}
while  It\^{o}'s formula yields
\begin{equation}\label{existence thm 2}
  \begin{split}
    &E\left[e^{r(t)}\|u(t)\|^2-e^{r(T)}\|u(T)\|^2
        \right]\\
    =&\
         E\biggl[\int_t^{T} e^{r(s)}\biggl( \inptv{2\Psi(s),u(s)}\\
    &\quad
     -\|Z(s)\|^2-\left(K+\frac{4}{\lambda}\|v(s)\|_V^2+K\rho^2(v(s))\right)
     \|u(s)\|^2\biggr)\,ds\biggr].
  \end{split}
\end{equation}
By substituting \eqref{existence thm 2} into \eqref{exsitence thm 1}, we
get
\begin{equation}\label{existence thm 3}
  \begin{split}
    E\bigg[\int_0^{T}\varphi(t)&\bigg(
    \int_t^{T} e^{r(s)}\bigg(
        2\inptv{\Psi-\Phi(s,v(s),Z(s)),u(s)-v(s)}\\
    &-\Big(K+\frac{4}{\lambda}\|v(s)\|_V^2+K\rho^2(v(s))\Big)\|u(s)-v(s)\|^2
        \bigg)ds\bigg)dt\bigg]\leq 0.
  \end{split}
\end{equation}
Take $v=u-\gamma\phi w$ for $\gamma>0$, $w\in V$ and $\phi \in L^{\infty}(\Omega\times [0,T],\mathcal {P},\mathbb{R})$. Then we divide by
$\gamma$ and let $\gamma\rightarrow 0$ to derive that
\begin{equation}
  E\bigg[\int_0^{T} \varphi(t)\Big(
    \int_t^{T} e^{r(s)}\phi(s)\Big(
        2\inptv{\Psi-\Phi(s,u(s),Z(s)),w}
        \Big)ds\Big)dt\bigg]\leq 0.
\end{equation}
By the arbitrariness of $\varphi$,$\phi$ and $w$, we have
$$\Gamma=\Psi=\Phi(\cdot,u,Z),~~a.e. \textrm{ on } \Omega\times [0,T]$$
In view of~\eqref{thm main prof step31}
  and keeping in mind the fact $\bar{u}=u$ $dt\times \mathbb{P}$-$a.e.$, we have
\begin{equation}
  u(t)=\xi+\int_t^T\Phi(s,u(s),Z(s))\, ds-\int_t^T Z(s)\, dW_s.
\end{equation}
Hence, by Remark \ref{rmk def} we conclude that $(u, Z)$ is a strong
solution to the 2D BSNSE problem \eqref{pb}.

\tbf{Step 4}. We shall prove the uniqueness. Suppose that there are two
strong solutions $(u(\cdot), Z(\cdot))$ and $(v(\cdot), Y(\cdot))$ to the
problem \eqref{pb} corresponding to the same terminal data $\xi$. Then
\begin{align*}
   u(t) - v(t) & = \int_t^T (- \nu Au(s) + \nu Av(s))\, ds + \int_t^T (B(u(s)) - B(v(s)))\, ds \\
   & +\int_t^T (JZ(s)-JY(s))\,ds+ \int_t^T (f(s, u(s), Z(s))  - f(s, v(s), Y(s)))\, ds \\
   &- \int_t^T (Z(s) - Y(s))\, dW_s, \; \; t \in [0, T], \; a.s..
\end{align*}
Define
$$
R(t)=R(\omega,t)=\int_0^t(K+\frac{4}{\lambda}\|v(\omega,s)\|_V^2+K\rho^2(v(\omega,s)))\,ds,
\quad (\omega,t)\in \Omega\times[0,T].
$$
Then in view of Lemma \ref{L:difest} and by It\^{o}'s formula (for instance, see \cite[Theorem 4.2.5]{PR07}), we have
\begin{equation}
  \begin{split}
    &E_{\mathscr{F}_{t}}\left[ e^{R(t)}\|u(t)-v(t)\|^2\right]\\
    =&E_{\mathscr{F}_{t}}\bigg[\int_t^{T} e^{R(s)}\Big(
    2\inptv{\Phi(s,u(s),Z(s))-\Phi(s,v(s),Y(s)),u(s)-v(s)} \\
     &
     -\|Z(s)-Y(s)\|^2 -\left(\kappa+\frac{4}{\lambda}\|v(s)\|_V^2+K\rho^2(v(s))\right)\|u(s)-v(s)\|^2\Big)\,ds\bigg]  \\
    \leq&
        E_{\mathscr{F}_{t}}\bigg[\int_t^{T} e^{R(s)}\Big(
        -\frac{\bar{\lambda}^2-1}{2\bar{\lambda}^2}\|Z(s)-Y(s)\|^2
        -\lambda\|u(s)-v(s)\|_V^2
        \Big)\,ds\bigg].
  \end{split}
\end{equation}
Thus,
\begin{equation*}
\begin{split}
E_{\mathscr{F}_t}\!\!\!\left[ e^{R(t)}\|u(t)-v(t)\|^2 +\int_t^{T}\!\!\!\! e^{R(s)}\Big(
        \frac{\bar{\lambda}^2-1}{2\bar{\lambda}^2}\|Z(s)\!-\!Y(s)\|^2
        \!+\!\lambda\|u(s)\!-\!v(s)\|_V^2
        \Big)ds \right] \leq 0,
\end{split}
\end{equation*}
which implies
$$\textrm{ for any } t\in [0,T],u(t)-v(t)=0\textrm{ in } H,~a.s. \textrm{ and }
E\left[\int_0^T\|Z(s)-Y(s)\|^2 ds\right]=0.$$ By the continuity of $u$ and
$v$, we have $\|(u-v,Z-Y)\|_{\ms{M}}=0,$ from which we  conclude that
$(u,Z)$ is only a modification of $(v,Y)$ in $(L_{\ms{F}}^2
   (\gw; C([0, T]; V)) \cap
    \ms{L}_{\ms{F}}^2 (0, T; D(A)))\times \ms{L}_{\ms{F}}^2 (0, T; V)$.
 We complete the proof.
\end{proof}


\end{document}